# Numerical analysis of pit-to-crack transition under corrosion fatigue using a stochastic pit generation algorithm

M. Mokhtari, X. Wang & J. Amdahl
*Centre for Autonomous Marine Operations and Systems (AMOS), Department of Marine Technology, Norwegian University of Science and Technology (NTNU), Norway*

ABSTRACT: Corrosion fatigue is a major threat to the integrity of marine structures. However, the combined damaging effect of corrosion and fatigue is not well understood due to the uncertainties associated with the stochastic nature of pitting corrosion. A pitting corrosion defect could have a quite complex morphology that varies randomly from one case to another. Nevertheless, in numerical corrosion fatigue studies, the complex pit morphology is often idealized using an overly simplified geometry with a smooth surface because of the difficulties involved in modelling defects with the irregular and random shapes of pitting corrosion. The present study investigates the influence of such geometrical simplifications on the results of numerical corrosion fatigue analyses. For this purpose, an isolated complex-shaped pit, generated using a hierarchical stochastic algorithm scripted in Python and linked with Abaqus/CAE, is developed in a Q235 steel dog-bone specimen. The numerical results obtained from this model are compared with those from another model containing an idealized counterpart of the irregular pit. A discussion on the effect of pit morphology on the stress/strain history and distribution and fatigue crack initiation is presented.

## 1 INTRODUCTION

Exposed to aggressive marine environments, marine structures and particularly their steel components are prone to corrosion. Although there are different corrosion protection methods, corrosion still occurs, accumulates, and results in accidents (Maureen et al. 2013). Both the economic and environmental losses due to corrosion are enormous (Gudze and Melchers 2008, Odusote et al. 2021, Pipeline Significant Incident 20 Year Trends | PHMSA, Zhang et al. 2018). Typically, pitting corrosion, a localized dissolution of metals usually initiated after the breakdown of protective coating or paint (Bhandari et al. 2015, Szklarska-Smialowska 1986), is the most common and concerning type of corrosion in the marine industry. Many research studies have investigated the effect of different parameters on pitting corrosion in marine environments (Gudze and Melchers 2008, Jeffrey and Melchers 2009, Melchers 2003a, b, 2009, 2018, 2019, Soares et al. 2009).

Within corrosion pits, there are often significant stress concentrations and large strain gradients. Therefore, under cyclic loading, a fatigue crack can be initiated from the pit site. The transition from pitting corrosion to fatigue crack has been studied for several decades by many researchers, yet it is not well understood. This is mostly due to the stochastic nature of pitting corrosion and its morphology. In early research on this topic, Kondos (1989) suggested a competition model between pitting corrosion and crack propagation. He simplified the problem by modelling the pit as an elliptical crack while the pitting corrosion can have quite complex 3-dimensional (3D) geometry. It was assumed that pit-to-crack transition occurs when the fatigue crack growth rate exceeds the pitting corrosion growth rate. Several studies followed this competition model (Dolley et al. 2000, Rokhlin et al. 1999, Zhou and Turnbull 1999). Further research considered the effect of pits, their sizes and occasionally their shapes on the corrosion fatigue to develop models for the pit-to-crack transition yet again considerable simplifications were made in modelling the pit and/or crack shape (Medved et al. 2004, Sadananda and Vasudevan 2020, Zhao et al. 2018). The stress concentration factor caused by pitting corrosion is analysed in several other studies and again simplistic geometries are employed to model the shape of the individual pits of interest (Cerit 2019, Huang et al. 2014, Shojai et al. 2022, Xu and Wang 2015).

Despite numerous investigations on pitting corrosion fatigue, the combined damaging effect of corrosion and fatigue is still not well understood. As noted earlier, this is mostly due to the uncertainties associated with the stochastic nature of pitting corrosion. The morphology of the pitting corrosion defect could be quite complex such that the pit topography could vary randomly from one case to another. It is demonstrated that the irregular pit morphology is often caused by the coalescence of smaller pits that have grown and merged into a bigger pit (Melchers 2015, 2018, Mokhtari and Melchers 2016, 2018, 2019, 2020). Secondary pits could emerge and grow inside an existing pit resulting in the so-called 'pits-within-pits' phenomenon, which could lead to a significantly rough surface inside an isolated pit. As discussed earlier, in most numerical studies, the irregular pit morphology is idealized using a simple geometry with a smooth surface due to the difficulties to model the irregular and random shapes of pitting corrosion (Cui et al. 2020, Guzmán-Tapia et al. 2021, Hu et al. 2021, Jiang et al. 2022, Jie and Susmel 2020, Moghaddam et al. 2019, Rezende et al. 2022, Xu et al. 2021). This is also the case in the artificial modelling of longer-term pits in both numerical and experimental studies (Abdalla Filho et al. 2014, Bhardwaj et al. 2022, Cerit 2013, Chiodo and Ruggieri 2009, Hu et al. 2016, Liu et al. 2021, Netto et al. 2005, Wu et al. 2022, Xu and Cheng 2013, Xu et al. 2017, Yeom et al. 2015, Zhang, Rathnayaka, et al. 2017, Zhang and Tian 2022, Zhang, Shi, et al. 2017, Zuñiga Tello et al. 2019).

With a new algorithm recently developed by Mokhtari and Melchers (2018, 2020), it is now possible to rapidly and realistically generate the complex morphology of pitting corrosion in both microscopic and macroscopic scales. Using this algorithm, many pitting corroded samples could be efficiently generated in research studies. This facilitates understanding the effect of pit morphology, including its irregularities and stochastic geometrical properties, on the structural and mechanical behaviour of corroded components. Using this stochastic pit generation algorithm, the present study investigates the effect of pit morphology on pit-to-crack transition with a focus on the consequences of pit shape idealization based on simple geometries in the numerical analyses.

## 2 METHODOLOGY

### 2.1 Finite element modelling

Three different finite element models of a dog-bone tensile test sample (Figure 1), prepared as per GB/T 3075-2008 (GB/T 3075-2008), are developed using Abaqus/CAE 2019 (Dassault Systèmes SIMULIA User Assistance, Abaqus 2019) and solved by the nonlinear implicit solver Abaqus/Standard. These models include one intact and two corroded specimens with isolated pitting corrosion defects. A realistic complex-shaped pitting corrosion profile is artificially generated in the first corroded specimen using a hierarchical stochastic algorithm implemented in several studies (Mokhtari and Melchers 2016, 2018, 2019, 2020). This complex-shaped defect was idealized in the second corroded model using an ellipsoidal profile to investigate the effect of pit shape idealization on the numerical results. In other words, the isolated pits have the same maximum depth, $d$, and width, $w$, in both of the corroded specimens but are modelled with two different shapes, ellipsoidal and irregular/complex-shaped. Note that in the irregular pit studied herein, $w$ is larger than the pit length, $l$. For pit shape idealizations, often the larger dimension (between $l$ and $w$) is used to define the pit diameter on the surface, and thus $w$ is opted as the ellipsoidal pit diameter on the specimen surface, $D$.

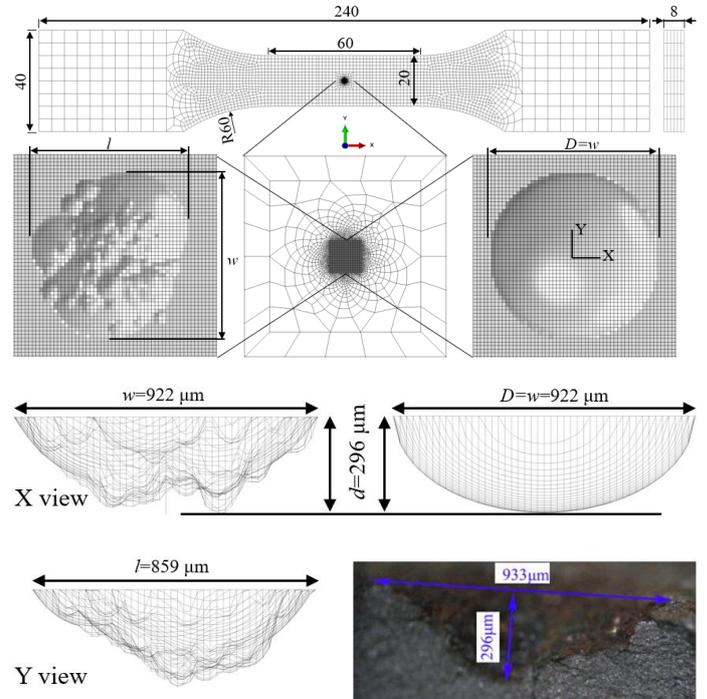

Figure 1. Finite element models with a complex-shaped pit (left) and its idealized counterpart (right). An example of a real critical pit in the dog-bone sample tested in (Xu and Wang 2015) is shown in the right-bottom corner. Dimensions in the dog-bone sample drawing are in mm.

The specimens are made of Q235 steel with the chemical composition given in Table 1. This type of steel is widely used in marine structures and steel bridges. The basic mechanical properties of the Q235 steel investigated in the present study (Table 2) are obtained from uniaxial tensile tests carried out in (Xu and Wang 2015). Given that the stress-strain curve of the Q235 steel is not provided in the reference study, a trilinear elastic-plastic material model is employed to define the material behaviour in the finite element models of this study. Note that the reference study also employed a trilinear elastic-plastic material model to numerically simulate the stress/strain distributions. Linear elastic behaviour is assumed with Young's modulus of 198 GPa and Poisson's ratio of 0.3 (Xu and Wang 2015). The intact finite element model is discretized into 12510 linear hexahedral elements with reduced integration (C3D8R). The corroded samples, however, are constructed with 54764 C3D8R elements. The element size and mesh configuration for each sample (Figure 1) resulted from a mesh sensitivity analysis with respect to the maximum von Mises stress.

Table 1. Q235 chemical composition (wt%) (Xu and Wang 2015).

| C | Si | Mn | P | S | Cr | Ni | Cu | N |
|---|---|---|---|---|---|---|---|---|
| 0.2 | 0.35 | 1.4 | 0.045 | 0.045 | 0.3 | 0.3 | 0.3 | 0.008 |

Table 2. Mechanical properties of Q235 obtained from uniaxial tensile testing at room temperature (Xu and Wang 2015).

| Yield strength $\sigma_y$ (MPa) | Ultimate strength $\sigma_u$ (MPa) | Elongation $\varepsilon_u$ (%) | Young's modulus $E$ (GPa) | Poisson's ratio $v$ |
|---|---|---|---|---|
| 312.3 | 458.0 | 38.2 | 198 | 0.3 |

The morphology of the artificial irregular pit in Figure 1, is produced automatically by a Python code linked with Abaqus/CAE that incrementally cuts spherical-cap-shaped geometries from the model surface resulting in an amalgamation of many smaller pits that are merged into a larger irregular pit. This cutting process is done in different hierarchical levels. In the first level, large pits are cut from the intact surface and in the other levels, secondary pits or sub-pits are cut from the corroded surface generated in the first level. The pit parameters (i.e., centre coordinates and radii) follow stochastic distributions. More information about the irregular pit generation algorithm and how to replicate the geometry of a real complex-shaped pit is provided in (Mokhtari and Melchers 2018, 2020).

All three specimens have fixed boundary conditions at one end of the dog-bone sample while the other end is loaded such that the maximum and minimum nominal stresses, $\sigma_{max}$ and $\sigma_{min}$, in the gauge section of the sample (i.e., the middle, narrowed section of the dog-bone specimen) are respectively 260.0 and 26.0 MPa. This is consistent with the reference study (Xu and Wang 2015).

## 2.2 Fatigue life analysis

After having the finite element analysis (FEA) completed, the numerical results (i.e, stress/strain data) are imported into FE-SAFE (3DS Simulia 2019b) for fatigue life analysis. The load history, surface roughness and material properties can be defined in FE-SAFE. To determine the fatigue life, different criteria based on stress, strain or energy can be used. For multiaxial fatigue, the critical plane method developed by Brown and Miller (Brown and Miller 1973) reliably simplifies the multiaxial complications into a uniaxial problem (Yang and Sun 2022). Further development of this method in (Kandil et al. 1982) led to the Brown-Miller combined strain criterion that is known to return the most realistic life estimates for ductile metals (3DS Simulia 2019a). According to the Brown-Miller criterion, maximum fatigue damage occurs on the maximum shear strain plane and the damage is a function of both the strain normal to this plane, $\varepsilon_n$, and the maximum shear strain, $\gamma_{max}$. Therefore, the conventional strain-life equation,

$$\frac{\Delta\varepsilon}{2} = \frac{\sigma'_f}{E}\left(2N_f\right)^b + \varepsilon'_f\left(2N_f\right)^c \qquad (1)$$

can be rewritten with the shear strain and normal strain amplitudes as

$$\frac{\Delta\gamma_{max}}{2} + \frac{\Delta\varepsilon_n}{2} = C_1\frac{\sigma'_f}{E}\left(2N_f\right)^b + C_2\varepsilon'_f\left(2N_f\right)^c \qquad (2)$$

where, $C_1$ and $C2$ are material constants and $\varepsilon$, $\varepsilon'_f$, $c$, $\sigma'_f$, $b$, and $N_f$ are respectively strain, fatigue ductility coefficient, fatigue ductility exponent, fatigue strength coefficient, fatigue strength exponent, and fatigue life. The first term on the right side of Equations 1 and 2 is the elastic component developed by Basquin (Basquin 1910), and the second term is the plastic component proposed by Coffin-Manson (Coffin Jr 1954, Manson 1953) in the strain-life approach (Prueter 2021). The $C_1$ and $C2$ constants can be found by solving Equation 2 for a uniaxial plane stress problem, where principal strains $\varepsilon_1$, $\varepsilon_2$, $\varepsilon_3$ are related to each other as follows.

$$\varepsilon_3 = \varepsilon_2 = -v\varepsilon_1 \qquad (3)$$

From Mohr's strain circle, $\gamma_{max}$ and $\varepsilon_n$, can be expressed as

$$\gamma_{max} = \frac{\varepsilon_1 - \varepsilon_3}{2} \qquad (4)$$

$$\varepsilon_n = \frac{\varepsilon_1 + \varepsilon_2}{2} \qquad (5)$$

Through Equations 2–5 and using the values of 0.3 and 0.5 for the elastic and plastic Poisson's ratios (i.e., $v_e$=0.3 and $v_p$=0.5), $C_1$ and $C_2$ are found to be 1.65 and 1.75, respectively. Therefore, the Brown-Miller combined strain criterion with Morrow mean stress correction (Morrow 1968) can be expressed as

$$\frac{\Delta\gamma_{max}}{2} + \frac{\Delta\varepsilon_n}{2} = 1.65\frac{(\sigma'_f - \sigma_{n,m})}{E}\left(2N_f\right)^b + 1.75\varepsilon'_f\left(2N_f\right)^c \qquad (6)$$

where $\sigma_{n,m}$ is the mean normal stress on the plane (3DS Simulia 2019a). Using Equation (6) (i.e., the Brown-Miller combined strain criterion with Morrow mean stress correction), multiaxial fatigue life analyses using the fatigue properties of Q235 given in Table 3 are performed in the present study.

Table 3. Q235 fatigue properties from (Fan et al. 2018).

| Cyclic strain hardening exponent $n'$ | Cyclic strain hardening coefficient $K'$ (MPa) | Fatigue ductility coefficient $\varepsilon'_f$ |
|---|---|---|
| 0.125 | 895.0 | 2.63 |
| Fatigue ductility exponent, $c$ | Fatigue strength coefficient, $\sigma'_f$ (MPa) | Fatigue strength exponent, $b$ |
| -0.89 | 1010.0 | -0.1113 |

# 3 RESULTS AND DISCUSSION

## 3.1 *Validation*

FE-SAFE calculates the fatigue life using the stress/strain data in each finite element. Therefore, if the element size is small enough (i.e., in the microscopic range), the fatigue life of the element could be taken as fatigue crack initiation life, $N_i$, assuming the crack has the same size as the elements with the minimum fatigue life. In the present study, the size of surface elements inside the pitting region is 20 μm on the surface and around 17–87 μm along the thickness while the element size at the gauge region of the intact model is 1.33 mm. This means that the fatigue life of an element in the intact model could be taken as the number of cycles to develop a crack of at least 1.33 mm which is around 17% of the specimen thickness. In addition, Figure 2 shows that the elements through the thickness of the intact model have similar fatigue lives, from 6.784 to 6.878 logarithmic cycles (i.e., 6.08–7.55 million cycles). As a result, the fatigue life produced by FE-SAFE for the intact model is deemed to be approximately the total life of the specimen, $N_t$. Given that the crack initiation life was not measured in the reference study (Xu and Wang 2015) and that FE-SAFE cannot estimate the total life of the corroded specimens studied herein (note total life estimation is out of the scope of this study), the total life of the intact model is used to validate the numerical model of the present study. Details of the fatigue life analysis of the intact model in both numerical and experimental studies are given in Table 4. The fatigue lives of the elements in the gauge section of the intact model are between 6.08 and 7.55 million cycles which is in good agreement with the experimental fatigue life of 6.73 million cycles.

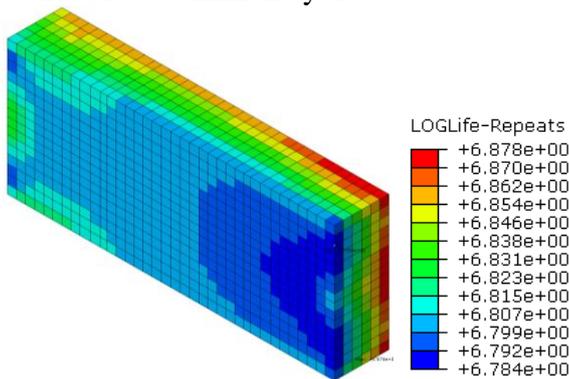

Figure 2. Logarithmic life calculated by FE-SAFE for the gauge section of the intact model.

Table 4. Fatigue analysis details of the intact model.

| Parameter | Experimental (Xu and Wang 2015) | Numerical |
| --- | --- | --- |
| Min cross-section area (mm$^2$) | 155.32 | 160.0 |
| Surface roughness (μm) | 16.868 | ≤16 |
| $\sigma_{max}$ (MPa) | 260.0 | 260.0 |
| Load ratio, $R$ | 0.1 | 0.1 |
| $N_t$ ($\times 10^6$) | 6.73 | 6.08–7.55 |

## 3.2 *Stress/strain history and distribution*

History plots for the equivalent plastic strain, $\bar{\varepsilon}_p$, von Mises stress, $\sigma_{von}$, and longitudinal stress, $\sigma_{xx}$, are plotted and discussed in this section. These plots will show the variations of stress and strain in the critical element for the first four cycles. The critical element here is opted to be the element with the maximum equivalent strain (Figure 3) for consistency with the literature and the reference study (Xu and Wang 2015). However, the location of maximum plastic strain may not necessarily be the crack initiation location or the only crack initiation location as will be shown in the next section.

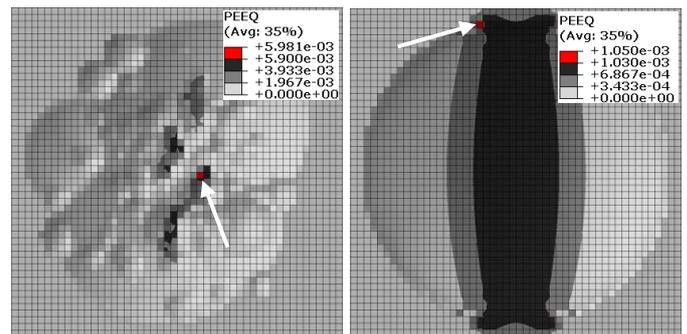

Figure 3. Equivalent plastic strain distribution and the location of the critical element shown with the arrow pointing to the element in red.

The significant effect of pit shape idealization on the equivalent plastic strain is evident in Figure 4 such that the irregular pit has caused 279% higher plastic strain at the pit site as opposed to the ellipsoidal pit model. This huge difference, of course, is due to the very rough base of the irregular pit studied in the present. Plastic strain in an irregular pit with a smoother base would evidently be smaller and closer to that caused by an ellipsoidal pit. Having said that, the element size on the pitting surface herein is 20 μm, meaning that surface roughness below this value that may cause even larger plastic strains could not be captured. After all, these results demonstrate the necessity for realistic modelling of pit morphology in the numerical analysis of irregular pits.

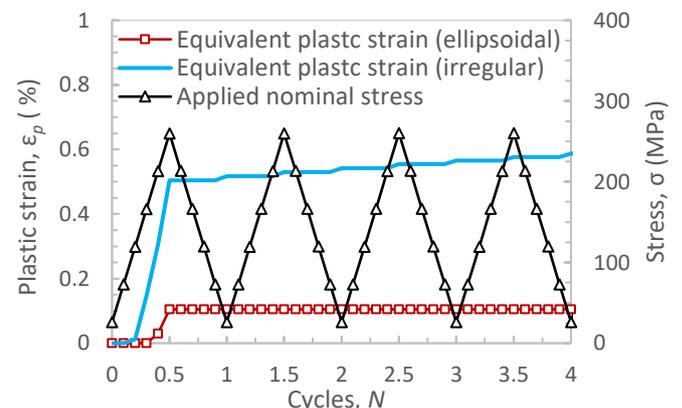

Figure 4. Variation of the equivalent plastic strain in the critical element during the first four cycles.

Figure 4 clarified the importance of realistic modelling of irregular pits geometry. However, since the common practice for the numerical fatigue life analysis of pitting corroded specimens usually uses the maximum stress obtained from a monotonic load, the stress history is also studied in this section. In Figure 5, although the peak von Mises stresses that occurred at $\sigma_{max}$ moment are almost the same for both irregular and ellipsoidal cases, their residual stresses that occurred at $\sigma_{min}$ moment are quite different. This can be seen better when the axial stresses are compared for the two pit modelling methods (Figure 6). Despite the tension-tension cyclic load applied to the specimen, the critical element experiences an almost fully reversed cyclic load due to the very large compressive residual stress inside the irregular pit. This is not the case for the ellipsoidal pit. Consequently, while the frequency of von Mises stress in the elliptical pit is the same as that of the nominal applied stress, $\sigma_{max}$, it is doubled when the pit has an irregular shape. In addition, at the beginning of the loading, before there is any plastic deformation in the ellipsoidal pit ($N<0.4$), the stress results obtained from the two models are very different. Therefore, stress-based analysis of a specimen under monotonic load and plastic deformation could be quite misleading as both irregular and ellipsoidal pits develop almost the same peak stress. As such, one may falsely conclude that due to the similar peak stress results, ellipsoidal idealization of an irregular pit shape is a reliable practice.

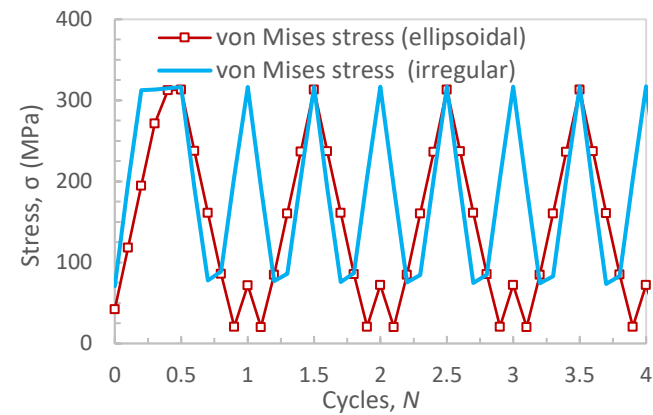

Figure 5. Variation of the von Mises stress in the critical element during the first four cycles.

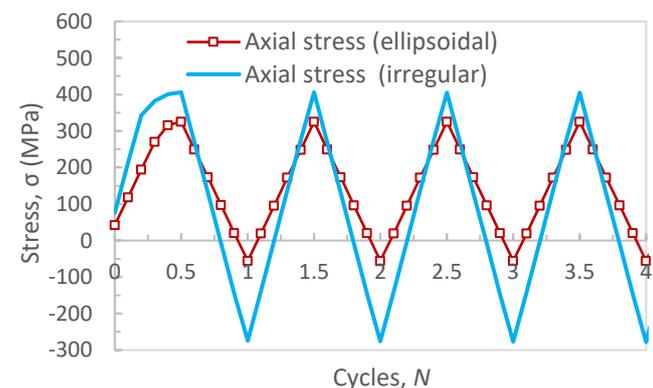

Figure 6. Variation of the axial stress in the critical element during the first four cycles.

## 3.3 *Fatigue life analysis*

Due to the considerable localized plastic strain in the corroded specimens, a nonlinear fatigue life analysis in FE-SAFE was carried out using both stress and strain data. According to Figures 7 and 8, the crack initiation life caused by the irregular pit ($N_i = 10^{4.151} = 14.2\times10^3$ cycles) is significantly smaller than that induced by its ellipsoidal counterpart ($N_i = 10^{5.412} = 258.2\times10^3$ cycles).

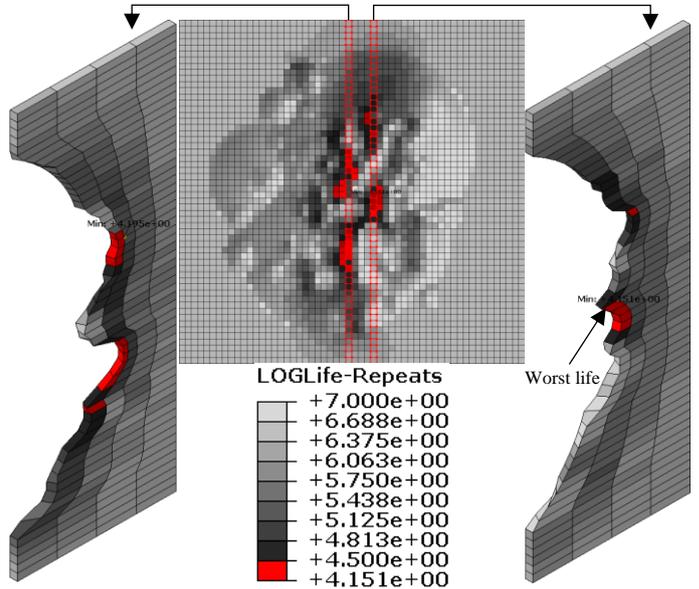

Figure 7. Logarithmic fatigue life distribution in the irregular pit. The red shows crack sites at $31.6\times10^3$ cycles (Worst life=$14.2\times10^3$ cycles).

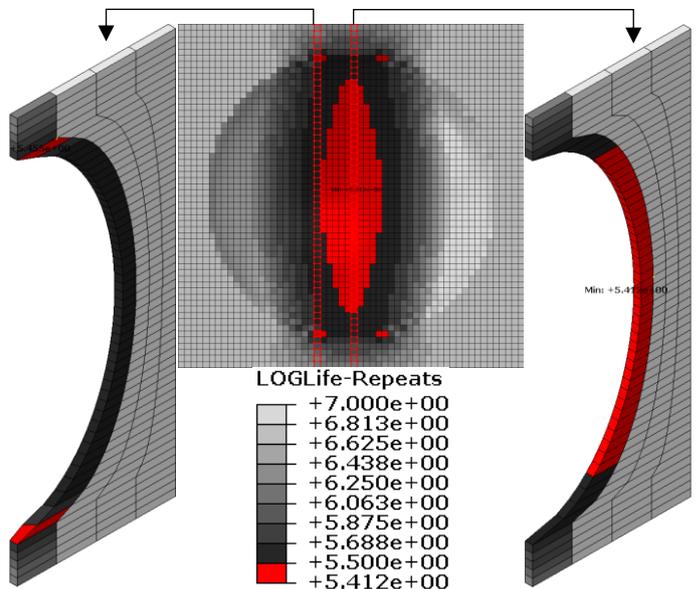

Figure 8. Logarithmic fatigue life distribution in the ellipsoidal pit. The red shows crack sites at $316.2\times10^3$ cycles (Worst life=$258.2\times10^3$ cycles).

In the experimental tests (Xu and Wang 2015), only the total fatigue life of the corroded specimen is reported ($N_t = 571.0\times10^3$), and there is no information on $N_i$. However, the small $N_i/N_t$ ratio ($N_i/N_t = 2.5\%$) predicted using the irregular pit shape is consistent with the literature (DuQuesnay et al. 2003, Gruenberg et al. 2004, Sankaran et al. 2001, Van der Walde et al. 2005, Van der Walde and Hillberry 2008) and the reference study (Xu and Wang 2015) reporting that very few cycles relative to the total life are needed to initiate a microcrack such that $N_i$ could be negligible in thick specimens. On the other hand, $N_i/N_t$ =45.2% predicted using the ellipsoidal pit shape is far too large and unrealistic.

The crack initiation sites in the irregular pit also differ from those in its ellipsoidal counterpart, which was expected. For the ellipsoidal pit, there are two different types of crack initiation sites, one located at the pit base where the pit reaches its maximum depth and another one positioned at the pit mouth where the elements with maximum equivalent plastic strain are situated. In contrast, there are several crack initiation sites in the complex-shaped pit, and all are located in the secondary pits at the base of the irregular pit. These secondary pits are not necessarily at the maximum pit depth region. In fact, the element with the worst life is located 62 μm higher than the maximum pit depth point near the mouth of a secondary pit that is not the deepest one. These results show that neither the site nor the life of fatigue crack initiation can be correctly predicted by models that use oversimplified pit shape idealization.

## 4 CONCLUSIONS

The effect of pit morphology on pit-to-crack transition under corrosion fatigue was investigated with a focus on the consequences of pit shape idealization in the numerical analyses. For this purpose, an isolated complex-shaped/irregular pit was developed in a finite element model of a Q235 steel dog-bone specimen. Then, the shape of this pit was idealized into an ellipsoidal geometry in another finite element model. The following conclusions are derived by comparing the stress and strain histories and fatigue crack initiation lives obtained from the two finite element models:

With a nominal maximum stress of 260 MPa and $R$=0.1, large residual stresses were found in both the irregular and ellipsoidal pits. However, the residual stress in the irregular pit was found to be far larger and much more localized than that in the ellipsoidal one (338% larger in terms of von Mises stress).

The irregular pit caused 279% higher equivalent plastic strain in the first cycle compared to its ellipsoidal counterpart. This difference in plastic strain gradually increased in the following cycles.

In the elastic range, before any plastic deformation occurs, the maximum von Mises stress in the irregular pit was significantly larger (up to 60%) than that in its ellipsoidal counterpart.

Despite all the significant differences between the results obtained from the models with irregular and ellipsoidal pits, their maximum local stresses reached almost the same value upon exceeding the yield stress. Therefore, under plastic deformations, a stress-based monotonic load analysis does not show the effect of pit shape idealization and could be quite misleading.

Neither the site nor the life of fatigue crack initiation was correctly predicted by the idealized model with an ellipsoidal pit.

## 5 REFERNCES